\documentclass{article}
\usepackage[american]{babel}
\usepackage{amsopn}
\usepackage{amsfonts}
\usepackage{hyperref}
\usepackage{amsmath,amssymb}
\usepackage{fancyhdr}
\numberwithin{equation}{section}

\def\xL{{\rm L}}
\def\xW{{\rm W}}
\def\xCn{{\rm C}}

\def\indic{{\rm {\large 1}\hspace{-2.3pt}{\large
l}}}

\def\pre{{\mathfrak{Re}}}
\def\xR{{\mathbb R}}
\def\E{{\mathbb E}}
\def\xN{{\mathbb N}}

\def\xL{{\rm L}}
\def\xW{{\rm W}}
\def\xCn{{\rm C}}
\def\xdif{{\rm d}}
\def\xLtwo{{{\rm L}^2}}
\def\xHone{{{\rm H}^1}}

\def\xLn{{\rm L}}

\def\xHn{{\rm H}}
\def\xtr{{\rm tr}}
\def\H1{{\xHone(\R^d)}}

\def\indic{{\rm {\large 1}\hspace{-2.3pt}{\large
l}}}

\def\xim{{\rm im}}

\newtheorem{thrm}{Theorem}[section]
\newtheorem{lmm}[thrm]{Lemma}

\newtheorem{prpstn}[thrm]{Proposition}

\newtheorem{rmrk}[thrm]{Remark}

\begin{document}
\title{Stochastic nonlinear Schr\"odinger equations driven by a fractional noise\\ Well posedness, large deviations and support}


\date{}

\maketitle

\begin{center}
{\large Eric Gautier}\vspace{0.3cm}

Cowles Foundation, Yale University,\\
30 Hillhouse Avenue, New Haven CT 06520, USA\\
e-mail: eric.gautier@yale.edu\\
http://pantheon.yale.edu/~eg298
\end{center}

\begin{abstract}We consider stochastic nonlinear Schr\"odinger
equations driven by an additive noise. The noise is fractional in
time with Hurst parameter $H$ in $(0,1)$. It is also colored in
space and the space correlation operator is assumed to be nuclear.
We study the local well-posedness of the equation. Under adequate
assumptions on the initial data, the space correlations of the noise
and for some saturated nonlinearities, we prove a sample path large
deviations principle and a support result. These results are stated
in a space of exploding paths which are H\"older continuous in time
until blow-up. We treat the case of Kerr nonlinearities when
$H>\frac12$.
\end{abstract}

 \vspace{0.3cm}


\newpage

\section{Introduction}\label{s1}
Nonlinear Schr\"odinger (NLS) equations are a generic model for the
propagation of the enveloppe of a wave packet in weakly nonlinear
and dispersive media, see \cite{SS}. They appear for example in
nonlinear optics, hydrodynamics, biology, field theory, crystals,
Bose-Einstein condensates, Fermi-Pasta-Ulam chains of atoms.
Sometimes random perturbations have to be considered. In optics,
noise accounts for the spontaneous emission noise due to amplifiers
placed along the fiber line in order to compensate for loss in the
fiber. In the context of crystals or of Fermi-Pasta-Ulam chains of
atoms, noise sometimes accounts for thermal effects. Noises
considered are often either complex additive noises or real
multiplicative noises. In physics, the Gaussian space-time white
noise is broadly considered. It has not been possible so far to give
a mathematical meaning to the solutions of such equations. Noises
considered in mathematics are colored in space. Note that in optics,
because the time variable corresponds to space and the space
variable to some retarded time, noises considered for well-posedness
are indeed
colored in time.\\
\indent We consider here the case of a fractional additive noise.
Fractional noises, introduced by Mandelbrodt, have several
applications in hydrology, finance and telecommunications. They are
extensions of the Gaussian white noise ($H=\frac12$) and for
$H\ne\frac12$ the noises are colored in time. Up to our knowledge,
these noises have not been considered in Physics for such models.
Again, in optics the new correlations could account for correlations
in space. However, we consider such noises to show that the results
of \cite{EG1} can be extended to more general Gaussian noises. We
specify these particular fractional noises
for computational issues since we then know precisely the kernel.\\
\indent The stochastic NLS equations could be written with the It\^o
notations
\begin{equation}
\label{e1} i \xdif u -(\Delta u + f(u))\xdif t = \xdif W^H,
\end{equation}
where $u$ is a complex valued function of time and space and $W^H$
is a fractional Wiener process. The fractional noise is formally its
time derivative. The parameter $H$ is called the Hurst parameter. It
belongs to $(0,1)$. The space variables belong to the whole space
$\xR^d$. The initial datum $u_0$ is a function of a particular
Sobolev space based on $\xLtwo$.\\
\indent We consider pathwise weak solutions in the sense used in the
analysis of PDEs. More precisely, we are interested in mild
solutions which are such that
\begin{equation}
\label{e2}u(t)=U(t)u_0-i\int_0^tU(t-s)f\left(u(s)\right)ds
-i\int_0^tU(t-s)dW^H(s),
\end{equation}
where $\left(U(t)\right)_{t\in\xR}$ is the Schr\"odinger linear
group on some Sobolev space  $\xHn^{s}$ generated by the
skew-adjoint unbounded operator $\left(-i\Delta,\xHn^{s+2}\right)$.\\
\indent Nonlinearities of the form $f(u)=\lambda |u|^{2\sigma}u$
where $\lambda=\pm 1$ are often considered for NLS equations, they
are called Kerr nonlinearities. In that case the space of energy
$\xHone$ is of particular interest. It is such that the Hamiltonian
is well defined. It is also a space where the blow-up phenomenon is
usually studied, indeed localized functions and thus the variance
may be defined. These nonlinearities are Lipschitz on the bounded
sets of $\xHone$ iff $d=1$. In higher dimensions, the Strichartz
inequalities, see \cite{SS} allow to treat these nonlinearities. For
the stochastic equations driven by a Gaussian noise which is white
in time, it is proved in \cite{dBD1} that the Cauchy problem is
locally well-posed in $\xHone$ for every $\sigma$ when $d=2$ and
only when $\sigma<\frac{2}{d-2}$ for $d\geq3$. Also, for such values
of $\sigma$ and $\lambda=-1$, defocusing case, the Cauchy problem is
globally well-posed. It is proved considering the mass and
Hamiltonian which are invariant quantities of the deterministic
equation. In the focusing case when $\lambda=1$, solutions may
blow-up in finite time when $\sigma\geq\frac{2}{d}$, critical and
supercritical nonlinearities. Note that in \cite{dBD2}, theoretical
results on the influence of a noise on the blow-up phenomenon have
been obtained. Large deviations and a support theorem for such
equations is given in \cite{EG1}. In \cite{EG2}, we prove LDP for a
noise of multiplicative type. In \cite{DG,EG1}, we apply our results
to the problem of error in soliton transmission by analyzing the
optimal control problem that governs the rate of the exponential
decay to zero with the noise intensity of the probability of large
deviation events. Note that, for the fractional noise of this
article, the computations would then be almost untractable because
of the extra correlations in time. Also, in \cite{EG4}, we apply the
uniform LDPs to the study the problem of the exit from a domain of
attraction for weakly damped equations. We use the strong Markov
property which does not hold for fractional noises. We also use LDPs
to obtain estimates on the small noise
asymptotic of the exit times in \cite{EG1,EG2}.\\
\indent A fractional Brownian motion (fBm) is a centered Gaussian
processes with stationary increments
\begin{equation*}
\E\left(\left|\beta^H(t)-\beta^H(s)\right|^2\right)=|t-s|^{2H},\quad
t,s>0.\end{equation*} A cylindrical fractional Wiener process on a
Hilbert space consists formally of independent fractional Brownian
motions (fBm) (with possibly different Hurst parameters) on each
coordinate of a complete orthonormal system. It does not have
trajectories in the Hilbert space. Also, only images by
Hilbert-Schmidt mappings are such that the laws of the marginals are
{\it bona fide} Radon measures. We cannot expect the stochastic
convolution to make sense removing the Hilbert-Schmidt assumption
since the group has no global smoothing properties in the Sobolev
spaces based on $\xLtwo$. It is an isometry on such spaces. Also,
since we work in $\xR^d$, the stochastic convolution would have to
be space wise translation invariant which is not compatible with the
fact that it should be a process with paths in a Sobolev space based
on $\xLtwo$. Thus we assume that the fractional Wiener processes in
our equation is a direct image
via a Hilbert-Schmidt operator of a cylindrical fractional Wiener process on $\xLtwo$.\\
\indent In this article, we consider the semi-group approach
developped in \cite{DPZ}. It is well suited for stochastic NLS
equations where we use properties of this group. It also allows to
define the stochastic integration in infinite dimensions with the
well studied integration with respect to the one parameter and one
dimensional fBm. More precisely we use the approach to the
stochastic calculus with respect to the fractional Brownian motion
developed in \cite{AMN} for general Voltera processes which is based
on the Malliavin calculus. In that
case, the stochastic integral is a Skohorod integral.\\
\indent We prove that the Cauchy problem is locally well-posed. We
then prove, for particular saturated nonlinearities {\it i.e.}
nonlinearities which are locally Lipschitz, a sample path large
deviation principle (LDP) for the small noise asymptotic and a
support theorem in a space of exploding paths which are
$H'-$H\"older continuous on time intervals before blow-up with
$0<H'<H$. Though the H\"older regularity holds for the stochastic
convolution it cannot be transfered easily since the group is an
isometry. We impose additional regulariy of the initial datum and
suitable assumptions on the correlations in space of the noise. In
the last section we treat the case of the Kerr nonlinearities for
$H>\frac{1}{2}$ but do
not impose conditions to obtain H\"older continuous paths.\\
\indent It is certainly much more involved to treat multiplicative
noises. For example, the stochastic convolution is now anticipating.
It is for the same reason that we do not investigate the global
existence. Indeed, the It\^o formula applied to the Hamiltonian and
mass to a certain power as in \cite{dBD1} gives rise to anticipating
stochastic integrals. These questions will be studied in future
works.

\section{Preliminaries}\label{s2}
The space of complex Lebesgue square integrable functions $\xLtwo$
with the inner product defined by
$(u,v)_{\xLtwo}=\pre\int_{\xR}u(x)\overline{v}(x)dx$ is a Hilbert
space. For $r$ positive, the Sobolev spaces $\xHn^r$ are the Hilbert
spaces of functions $f$ of $\xLtwo$ such that their Fourier
transform $\hat{f}$ satisfy
$\int_{\xR^d}\left(1+|\xi|^2\right)^r|\hat{f}(\xi)|^2d\xi<\infty$.
If $I$ is an interval of $\xR$, $(E,\|\cdot\|_E)$ a Banach space and
$r$ belongs to $[1,\infty]$, then $\xL^{r}(I;E)$ is the space of
strongly Lebesgue measurable functions $f$ from $I$ into $E$ such
that $t\rightarrow \|f(t)\|_E$ is in $\xL^{r}(I)$. The integral is
the Bochner integral. The space of bounded operators from $B$ to
$C$, two Banach spaces, is denoted by $\mathcal{L}_c(B,C)$. The
space of Hilbert-Schmidt operators $\Phi$ from $E$ to $F$, two
Hilbert spaces, is denoted by $\mathcal{L}_2(E,F)$. It is a Hilbert
space when endowed with the norm
$\|\Phi\|_{\mathcal{L}_2(E,F)}^2=\xtr\Phi\Phi^*=\sum_{j\in\xN}\|\Phi
e_j\|_{F}^2$ where $\left(e_j\right)_{j\in\xN}$ is a complete
orthonormal system of $E$. We denote by $\mathcal{L}_2^{0,r}$ the
above space when $E=\xLtwo$ and
$F=\xHn^r$.\\
\indent When $A$ and $B$ are two Banach spaces, $A\cap B$ with the
norm defined as the maximum of the norms in $A$ and in $B$, is a
Banach space. A pair $(r,p)$ of positive numbers is called an
admissible pair if $p$ satisfies $2\leq p<\frac{2d}{d-2}$ when $d>2$
($2\leq p<+\infty$ when $d=2$ and $2\leq p\leq+\infty$ when $d=1$)
and $r$ is such that
$\frac{2}{r}=d\left(\frac{1}{2}-\frac{1}{p}\right)$. Given an
admissible pair $(r(p),p)$ and $T$ positive, the space
\begin{equation*}X^{(T,p)}=\xCn\left([0,T];\xHone\right)\cap
\xLn^{r(p)}\left(0,T;\xW^{1,p}\right),\end{equation*} is the space
considered to prove the local existence of solutions to the NLS equation with a Kerr nonlinearity.\\
\indent Also, we denote by $\xCn^{H'}\left([0,T];E\right)$ the space
of $H'-$H\"older $E$-valued continuous functions on $[0,T]$ embedded
with the norm
\begin{equation*}
\left\|f\right\|_{H',T}=\sup_{t\in[0,T]}\|f(t)\|_{E}+\sup_{t,s\in[0,T],t\ne
s}\frac{\|f(t)-f(s)\|_{E}}{|t-s|^{H'}}
\end{equation*}
where $E$ is a Banach space. The space
$\xCn^{H',0}\left([0,T];E\right)$ is the separable subset of the
above such that
\begin{equation*}
\lim_{|t-s|\rightarrow 0}\frac{\|f(t)-f(s)\|_{E}}{|t-s|^{H'}}=0.
\end{equation*}
\indent We denote by $x\wedge y$ the minimum of $x$ and $y$. A rate
function $I$ is a lower semicontinuous function. It is good if for
every $c$ positive, $\left\{x:\ I(x)\leq c\right\}$ is
compact.\vspace{0.3cm}

Volterra processes, see for example \cite{Deu}, are defined for $T$
positive as
\begin{equation*}
X(t)=\int_0^tK(t,s)d\beta(s),\quad
K\in\xLtwo\left([0,T]\times[0,T]\right),\ T>0,\ K(t,s)=0\ {\rm if}\
s>t.
\end{equation*}
The covariance of such a process is
\begin{equation*}
R(t,s)=\int_0^{t\wedge s}K(t,r)K(s,r)dr.
\end{equation*}
The covariance operator, when we consider the $\xLtwo(0,T)-$random
variables, has finite trace. It could be defined through the kernel
$R(t,s)$, {\it i.e.} for $h$ in $\xLtwo(0,T)$,
$Rh(t)=\int_0^TR(t,s)h(s)ds$. Also $R$ is such that $R=KK^*$ where
$K$ is the Hilbert-Schmidt operator defined for $h\in\xLtwo(0,T)$ by
$Kh(t)=\int_0^TK(t,s)h(s)ds=\int_0^tK(t,s)h(s)ds$ and $K^*$ is its
adjoint. These processes admit modifications with continuous sample
paths; they are Gaussian processes. $\xim\ R^{\frac{1}{2}}$, the
range of $R^{\frac{1}{2}}$ with the norm of the image structure, is
the reproducing kernel Hilbert space (RKHS) of the Gaussian measure
which is the law of the process on $\xLtwo(0,T)$. It is classical
that is equal to $\xim\ K$. It is also the RKHS of the measure on
$\xCn([0,T])$ since the restriction of the measure is again a
Gaussian probability measure and $\xCn([0,T])$ is a Banach space
continuously embedded in $\xLtwo(0,T)$. Also, it is known that the
RKHS of the measure on $\xCn([0,T])$ is isometric to the closure in
$\xLtwo(\mu)$, where $\mu$ is the Gaussian measure, of the dual of
the Banach space defined by means of the evaluation at points $t$
in $[0,T]$ and thus to the first Wiener chaos.\\
\indent It is now standard fact, see for example \cite{BC,DS1}, that
the rate function of a LDP for the family of Gaussian measures
defined as direct images of $\mu$ via the mapping
$x\mapsto\sqrt{\epsilon}x$ is given by $\frac12\|\cdot\|_{\xim\
R^{\frac12}}^2$ and that the support of the law of the Gaussian
measure is the closure of the RKHS for the norm of the Banach space.
We aim to transport such results to the law of the mild solution of
the stochastic NLS equation driven by such noises. We know that it
is more convenient to primarily prove large deviations for a
modification of the infinite dimensional stochastic convolution with
smooth sample paths.\vspace{0.3cm}

Defining the stochastic integration in Hilbert spaces only requires
to define the stochastic integration in dimension one. Such Volterra
processes are seldom martingales. However, since they are Gaussian
processes, the Skohorod integral can be defined. Let us recall, for
the sake of completeness, some aspects of integration with respect
to Volterra processes with the Malliavin calculus, see \cite{AMN}
for more details. Let us first consider another space well suited
for the Malliavin calculus that is also a RKHS. It may be seen as
generated by step functions on $[0,T]$; the stochastic integral of a
step function $\indic_{[0,T]}$ should coincide with the evaluation
at point $t$. The set of step functions is denoted by $\mathcal{E}$.
We consider the inner product defined by
\begin{equation*}
R(t,s)=\left<\indic_{[0,t]},\indic_{[0,s]}\right>_{\mathcal{H}}
=\left(K(t,\cdot)\indic_{[0,t]},K(s,\cdot)\indic_{[0,s]}\right)_{\xLtwo(0,T)}.
\end{equation*}
The linear operator $K_T^*$ from $\mathcal{E}$ into $\xLtwo(0,T)$ is
defined for $\varphi$ in $\mathcal{E}$ by
\begin{equation}\label{e4}
\left(K_T^*\varphi\right)(s)=\varphi(s)K(T,s)+\int_s^T(\varphi(t)-\varphi(s))K(dt,s).
\end{equation}
It is such that, for any $\varphi$ in $\mathcal{E}$ and $h$ in
$\xLtwo(0,T)$, we have
\begin{equation*}
\int_0^T\left(K_T^*\varphi\right)(t)h(t)dt=\int_0^T\varphi(t)(Kh)(dt).
\end{equation*}
The RKHS $\mathcal{H}$ is now obtained as the closure of
$\mathcal{E}$ with respect to the norm
$\|\varphi\|_{\mathcal{H}}=\|K_T^*\varphi\|_{\xLtwo(0,T)}$. The
operator $K_T^*$ is then an isometry between $\mathcal{H}$ and a
closed subspace of $\xLtwo(0,T)$; we represent $\mathcal{H}$ as
$\mathcal{H}=\left(K_T^*\right)^{-1}\left(\xLtwo(0,T)\right)$. The
above duality relation allows to extend integration with respect to
$Kh(dt)$ to integrands in $\mathcal{H}$. It also allows to define a
stochastic integration with respect to these Voltera processes for
integrands $\varphi$ in $\mathcal{H}$ as the Skohorod integral
\begin{equation*}
\delta^X(\varphi)=\int_0^T\left(K_T^*\varphi\right)(t)\delta\beta(t).
\end{equation*}
For deterministic integrands, it is the case for the stochastic
convolution with an additive noise, the integral could be written as
a It\^o integral
\begin{equation*}
\delta^X(\varphi)=\int_0^T\left(K_T^*\varphi\right)(t)d\beta(t).
\end{equation*}
\indent From now on we restrict our attention to the particular case
of the fBm. Enlarging if necessary the probability space the fBm may
be defined in terms of a standard Brownian motion
$\left(\beta(t)\right)_{t\geq0}$ via the square integrable
triangular kernel $K^H$, {\it i.e.} $K^H(t,s)=0$ if $s>t$,
\begin{equation*}
\beta^H(t)=\int_0^tK^H(t,s)d\beta(s),
\end{equation*}
where
\begin{equation}\label{e3}
K^H(t,s)=c_H(t-s)^{H-\frac{1}{2}}+c_H\left(\frac{1}{2}-H\right)\int_s^t(u-s)^{H-\frac{3}{2}}
\left(1-\left(\frac{s}{u}\right)^{\frac{1}{2}-H}\right)du,
\end{equation}
and
\begin{equation*}
c_H=\left(\frac{2H\Gamma\left(\frac{3}{2}-H\right)}{\Gamma\left(H+\frac{1}{2}\right)
\Gamma\left(2-2H\right)}\right)^{\frac{1}{2}}.
\end{equation*}
\eqref{e3} implies that
\begin{equation}\label{e3b}
\frac{\partial K^H}{\partial
t}(t,s)=c_H\left(\frac{1}{2}-H\right)(t-s)^{H-\frac{3}{2}}\left(\frac{s}{t}\right)^{\frac{1}{2}-H}.
\end{equation}
We now denote the kernel and the operator by $K$ instead of $K^H$
for the fBm.\\
\indent Also, we recall the following properties. The fBm has a
modification with $H'-$H\"older continuous sample paths where
$0<H'<H$; see for example \cite{DU}. Its covariance is given by
\begin{equation*}
\E\left(\beta^H(t)\beta^H(s)\right)=\frac{1}{2}\left(s^{2H}+t^{2H}-|s-t|^{2H}\right).
\end{equation*}
The increments are independent if and only if $H=\frac12$. The
covariance of future and past increments is negative if $H<\frac12$
and positive if $H>\frac12$. Thus we say that the fBm presents long
range dependence for $H>\frac12$ as the covariance between
increments at a distance $u$ decays as $u^{2H-2}$. Finally note that
these processes are also self-similar, {\it i.e.} the law of the
paths $t\mapsto\beta^H(at)$ where $a$ is positive are that of
$t\mapsto a^H\beta^H(t)$. The solution of the NLS equation also
display a self similar behavior, but in the space variable, near
blow-up for supercritical nonlinearities, see \cite{SS}.\\
\indent In the particular case of the fBm the stochastic integration
may also be defined by means of the fractional calculus. A rough
paths approach may also be considered, see for example \cite{CQ}. We
expect that this latter approach could allow to treat noises of
multiplicative type.
\begin{rmrk}\label{rk1}
Suppose that for multiplicative noises we were able to prove the
continuity of the solution with respect to the driving process at
the level of the rough paths, see \cite{Ly} for certain SDEs, then
LDP and support theorems follow from a contraction principle stating
large deviations for the rough paths of the driving process. It is
done in \cite{LeQiZh} for a SDE driven by the Brownian motion. LDP
for the rough paths of the fBm and for a Banach space valued Wiener
process are proved in \cite{MSS} and \cite{LLQ}. Also a rough paths
approach to a linear SPDE with analytical semigroup and for "smooth"
rough paths is given in \cite{LGT}.
\end{rmrk}
\indent We use several times the following property, that we may
check using \eqref{e4} and \eqref{e3}, that for $0<t<T$,
\begin{equation}\label{e5}
\left(K_T^*\indic_{[0,t]}\varphi\right)(s)=\left(K_t^*\varphi\right)(s)\indic_{[0,t]}(s).
\end{equation}
For smooth kernels such that $H>\frac12$, relation \eqref{e4} has
the simpler form
\begin{equation}
\left(K_T^*\varphi\right)(s)=\int_s^T\varphi(r)K(dr,s).
\end{equation}
The formulation in \eqref{e4} however allows to extend this
definition to singular kernels, {\it i.e.} when $H<\frac12$. For $H$
such that $H>\frac12$, the inner product in $\mathcal{H}$ of
$\varphi$ and $\psi$ is given by
\begin{equation*}
\begin{array}{rl}
\left<\varphi,\psi\right>_{\mathcal{H}}=&\int_0^T\int_0^T\varphi(u)\psi(v)\int_0^{u\wedge v}\frac{\partial K}{\partial u}(u,s)\frac{\partial K}{\partial v}(v,s)dsdudv\\
=&c_H^2\left(H-\frac{1}{2}\right)^2B\left(2-2H,H-\frac{1}{2}\right)\int_0^T\int_0^T\varphi(u)\psi(v)|u-v|^{2H-2}dudv,
\end{array}
\end{equation*}
from a computation given in \cite{AN}; $B$ denotes the Beta function.
It corresponds to the covariance of the stochastic integrals with respect to the fBm
\begin{equation*}
\E\left[\int_0^T\varphi(u)d\beta^H(u)\int_0^T\psi(v)d\beta^H(v)\right];
\end{equation*}
the space $\mathcal{H}$ is thus what would be a RKHS at the level of the noise in $\xLtwo(0,T)$ which covariance
is $\int_0^{u\wedge v}\frac{\partial K}{\partial u}(u,s)\frac{\partial K}{\partial v}(v,s)ds$.\\
\indent In Hilbert spaces, we assume that $W^H$ is the direct image
by a Hilbert-Schmidt operator $\Phi$ of a cylindrical fractional
Wiener process on $\xLtwo$, {\it i.e.} $W^H=\Phi W_c^H$. A
cylindrical fractional Wiener process on a Hilbert space $E$ is such
that for every orthonormal basis $\left(e_j\right)_{j\in\xN}$ of
$\xLtwo$ there exists independent fractional Brownian motions (fBm)
$\left(\beta_j^H(t)\right)_{t\geq0}$ such that
$W_c(t)=\sum_{j\in\xN}\beta_j^H(t)e_j$.\\
\indent Stochastic integration with respect to fractional Wiener
processes in a Hilbert space $F$, see for example \cite{TTV}, when
integrands are deterministic is defined as above but for step
functions multiplied by elements of the Hilbert space. It is such
that a scalar product by an element of $F$ is the one dimensional
stochastic integral of the scalar product of the integrand.
Operators $K_T^*$ are still well defined when the RKHS $\mathcal{H}$
is made of functions with values in $F$. Integrals of deterministic
bounded operator valued integrands $\Lambda$ from the Hilbert space
$E$ to $F$ are defined for $t$ positive as
\begin{equation*}
\int_0^t\Lambda(s)dW^H(s)=\sum_{j\in\xN}\int_0^t\Lambda(s)\Phi e_jd\beta_j^H(s)
=\sum_{j\in\xN}\int_0^t\left(K_t^*\Lambda(\cdot)\Phi e_j\right)(s)d\beta_j(s),
\end{equation*}
when $\left(\Lambda(t)\right)_{t\in[0,T]}$ is such that
\begin{equation*}
\sum_{j\in\xN}\int_0^T\|\left(K_T^*\Lambda(\cdot)\Phi
e_j\right)(t)\|_{F}^2dt<\infty.
\end{equation*}
Note that the duality relation \eqref{e2} still holds, the integral
is a Bochner integral, and that $K_T^*$ commutes with the scalar
product with an element of $F$. We now assume that $E=\xLtwo$ and
that $\left(e_j\right)_{j\in\xN}$ is a complete orthonormal system.
We may also check from \eqref{e2} that the linear group
$\left(U(t)\right)_{t\in\xR}$  on a Sobolev space based on $\xLtwo$
commutes with $K_T^*$.

\section{The stochastic convolution}\label{s3}
In this section we present a few properties of the stochastic
convolution.\\
\indent When we consider particular saturated nonlinearities, the
precise assumptions are given in the next section, we treat singular
kernels and state our results in spaces of H\"older continuous
functions. We thus make the following assumption\vspace{0.2cm}

\noindent {\bf Assumption (N1)}
\begin{equation*}
\Phi\ {\it belongs\ to\ }
\mathcal{L}_2\left(\xLtwo,\xHn^{1+2(H+\alpha)}\right)\ {\it with}\
\left(\frac12-H\right)\indic_{H<\frac12}<\alpha<(1-H)\indic_{H<\frac12}+\indic_{H\geq\frac12}.
\end{equation*}

\noindent This assumption is used along with the fact that for
$\gamma$ in $[0,1)$ and $t$ positive
\begin{equation}\label{PG}
\|U(t)-I\|_{\mathcal{L}_c(\xHn^{1+2\gamma},\xHn^1)}\leq2^{1-\gamma}|t|^{\gamma};
\end{equation}
it could be proved using the Fourier transform.\\
\indent When we consider Kerr nonlinearities when the space
dimension is such that $d>2$ we impose\vspace{0.2cm}

\noindent {\bf Assumption (N2)}
\begin{equation*}
\Phi\in\mathcal{L}_2^{0,2}\ {\it and}\ H>\frac12.
\end{equation*}
In \cite{dBD1}, the authors impose weaker assumptions on $\Phi$,
namely $\Phi\in\mathcal{L}_2^{0,1}$, and check the required
integrabilty of the stochastic convolution. It is more intricate for
a fractional noise. This integrability follows from the Strichartz
inequalities under (N2), however this assumption is certainly too
strong.\vspace{0.2cm}

\noindent Under (N1), the following result on the stochastic
convolution holds.
\begin{lmm}\label{l1}
The stochastic convolution $Z:\ t\mapsto\int_0^tU(t-s)dW^H(s)$ is
well defined. It has a modification in $C_{\infty}^{H,0}$ and
defines a $C_{\infty}^{H,0}-$ random variable. Moreover, the direct
images $\mu^{Z,T,H'}$ of its law $\mu^Z$ by the restriction on
$C_T^{H',0}$ for $T$ positive and $0<H'<H$ are centered Gaussian
measures.
\end{lmm}
\noindent {\bf Proof.} The stochastic convolution is well defined
since for $t$ positive
\begin{equation*}
\begin{array}{l}
\sum_{j\in\xN}\int_0^t\left\|\left(K_t^*U(t-\cdot)\Phi e_j\right)(u)\right\|_{\xHn^{1+2H}}^2du\\
=\sum_{j\in\xN}\int_0^t\left\|U(-u)\Phi e_jK(t,u)+\int_u^t\left(U(-r)-U(-u)\right)\Phi e_j K(dr,u)\right\|_{\xHn^{1+2H}}^2du\\
\leq2\|\Phi\|_{\mathcal{L}_2^{0,r+2H}}^2\int_0^tK(t,u)^2du+2\sum_{j\in\xN}\int_0^t
\left\|\int_u^t\left(U(-r)-U(-u)\right)\Phi e_j K(dr,u)\right\|_{\xHn^{1+2H}}^2du\\
\leq 2(T_1+T_2).
\end{array}
\end{equation*}
Note that we used the continuous embedding of $\xHn^{1+2(H+\alpha)}$
into $\xHn^{1+2H}$. The integral in $T_1$ is equal to
$\E[(\beta^H(t))^2]=t^{2H}$. Using \eqref{PG}, we obtain
\begin{equation*}
T_2\leq4^{1-\alpha}\|\Phi\|_{\mathcal{L}_2^{0,r+2(H+\alpha)}}^2c_H^2\left(\frac12-H\right)^2
\int_0^t\left(\int_u^t\left(r-u\right)^{H-\frac32+\alpha}\left(\frac{r}{u}\right)^{H-\frac12}dr\right)^2du
\end{equation*}
thus
\begin{equation*}
T_2\leq4^{1-\alpha}\|\Phi\|_{\mathcal{L}_2^{0,r+2(H+\alpha)}}^2c_H^2\left(\frac12-H\right)^2
\int_0^t\left(\int_u^t\left(r-u\right)^{H-\frac32+\alpha}dr\right)^2du,
\end{equation*}
the integral is well defined since $H-\frac32+\alpha>-1$. We finally
obtain
\begin{equation*}
T_2\leq\frac{4^{\frac12-\alpha}\|\Phi\|_{\mathcal{L}_2^{0,r+2(H+\alpha)}}^2}{H+\alpha}
\left(\frac{c_H\left(H-\frac12\right)}{H-\frac12+\alpha}\right)^2t^{2H+2\alpha}.
\end{equation*}
Note that when $H>\frac12$, the assumption on $\alpha$ is not
necessary, indeed the kernel is null on the diagonal and its
derivative is integrable. We could obtain directly
\begin{equation*}
T_2\leq\|\Phi\|_{\mathcal{L}_2^{0,r+2(H+\alpha)}}^2\int_0^tK(t,u)^2du=\|\Phi\|_{\mathcal{L}_2^{0,r+2(H+\alpha)}}^2t^{2H}.
\end{equation*}
\indent We now prove that for any positive $T$ and $0<H'<H$, $Z$ has
a modification in $C_T^{H',0}$. We prove that it has a modification
which is in $\xCn^{H''}\left([0,T],\xHn^1\right)$ for some $H''$
such that $H'<H''<H$. It will thus belong to
$\xCn^{H',0}\left([0,T],\xHn^1\right)$. Note that, as we are dealing
with a centered Gaussian process, upper bounds on higher moments
could be deduced from an upper bound on the second order moment; see
\cite{DPZ} for a proof in the infinite dimensional setting. It is
therefore enough to show that there exists positive $C$ and
$\gamma>H''$ such that for every $(t,s)\in[0,T]^2$.
\begin{equation*}
\E\left[\|Z(t)-Z(s)\|_{\xHn^{1}}^2\right]\leq C|t-s|^{2\gamma},
\end{equation*}
and then conclude with the Kolmogorov criterion.\\
\indent When $0<s<t$, we have
\begin{equation*}\begin{array}{rl}
Z(t)-Z(s)=&U(s)\left(U(t-s)-I\right)\sum_{j\in\xN}\int_0^T\left(K_T^*\indic_{[0,t]}(\cdot)U(-\cdot)\Phi e_j\right)(w)d\beta_j(w)\\
&+U(s)\sum_{j\in\xN}\int_0^T\left(\left(K_t^*U(-\cdot)\Phi e_j\right)(w)-\left(K_s^*U(-\cdot)\Phi e_j\right)(w)\right)d\beta_j(w)\\
=&\tilde{T}_1(t,s)+\tilde{T}_2(t,s).
\end{array}\end{equation*}
We have
\begin{equation*}
\begin{array}{l}
\E\left[\left\|\tilde{T}_1(t,s)\right\|_{\xHn^1}^2\right]\\
\leq\left\|U(t-s)-I\right\|_{\mathcal{L}_c\left(\xHn^{1+2(H+\alpha)},\xHn^1\right)}^2
\sum_{j\in\xN}\int_0^T\left\|\left(K_T^*\indic_{[0,t]}(\cdot)U(-\cdot)\Phi e_j\right)^2(w)\right\|_{\xHn^{1+2(H+\alpha)}}^2dw\\
\leq
C(T,H,\alpha)\|\Phi\|_{\mathcal{L}_2^{0,r+2(H+\alpha)}}^2|t-s|^{2(H+\alpha)},
\end{array}\end{equation*}
where $C(T,H,\alpha)$ is a constant, and
\begin{equation*}
\begin{array}{l}
\E\left[\left\|\tilde{T}_{2}(t,s)\right\|_{\xHn^1}^2\right]\\
=\sum_{j\in\xN}\int_0^T\|U(-u)\Phi e_jK(t,u)+\int_u^t\left(U(-r)-U(-u)\right)\Phi e_j K(dr,u)\\
\hspace{2.2cm} -U(-u)\Phi e_jK(s,u)-\int_u^s\left(U(-r)-U(-u)\right)\Phi e_j K(dr,u)\quad\ \|_{\xHn^{1}}^2du\\
\leq\sum_{j\in\xN}\left(\tilde{T}_{21}^j+\tilde{T}_{22}^j+\tilde{T}_{23}^j\right),
\end{array}
\end{equation*}
where, using the fact that the kernel is triangular,
\begin{equation*}
\begin{array}{rl}
\tilde{T}_{21}^j&=\int_0^s\left\|U(-u)\Phi e_j\left(K(t,u)-K(s,u)\right)+\int_s^t\left(U(-r)-U(-u)\right)\Phi e_j K(dr,u)\right\|_{\xHn^{1}}^2du,\\
\tilde{T}_{22}^j&=2\int_s^t\left\|U(-u)\Phi e_jK(t,u)\right\|_{\xHn^{1}}^2du\\
\tilde{T}_{23}^j&=2\int_s^t\left\|\int_u^t\left(U(-r)-U(-u)\right)\Phi
e_j K(dr,u)\right\|_{\xHn^{1}}^2du.
\end{array}
\end{equation*}
We have
\begin{equation*}
\begin{array}{rl}
\tilde{T}_{21}^j&=\int_0^s\left\|\int_s^tU(-r)\Phi e_j K(dr,u)\right\|_{\xHn^{1}}^2du\\
&=\left\|\Phi e_j\right\|_{\xHn^{s+\gamma}}^2\int_0^s\left(\int_s^t |K(dr,u)|\right)^2du\\
&=\left\|\Phi
e_j\right\|_{\xHn^{s+\gamma}}^2\int_0^s\left(K(t,u)-K(s,u)\right)^2du
\end{array}
\end{equation*}
thus
\begin{equation*}
\begin{array}{rl}
\tilde{T}_{21}^j&\leq\left\|\Phi e_j\right\|_{\xHn^{1}}^2\int_0^t\left(K(t,u)-K(s,u)\right)^2du\\
&\leq\left\|\Phi e_j\right\|_{\xHn^{1}}^2\E\left[\left(\beta^H(t)-\beta^H(s)\right)^2\right]\\
&\leq\left\|\Phi e_j\right\|_{\xHn^{1}}^2|t-s|^{2H},
\end{array}
\end{equation*}
and
\begin{equation*}
\begin{array}{rl}
\tilde{T}_{22}^j&=2\left\|\Phi e_j \right\|_{\xHn^{s+\gamma}}^2\int_s^tK(t,u)^2du\\
&=2\left\|\Phi e_j
\right\|_{\xHn^{s+\gamma}}^2\int_s^t\left(K(t,u)-K(s,u)\right)^2du
\end{array}
\end{equation*}
thus
\begin{equation*}
\begin{array}{rl}
\tilde{T}_{22}^j&\leq2\left\|\Phi e_j \right\|_{\xHn^{1}}^2\int_0^t\left(K(t,u)-K(s,u\right)^2du\\
&\leq2\left\|\Phi e_j \right\|_{\xHn^{1}}^2|t-s|^{2H},
\end{array}
\end{equation*}
finally the same computations as above shows hat when
$H-\frac{3}{2}+\alpha>-1$ (used for integrability issue when
$H<\frac12$), we have
\begin{equation*}
\begin{array}{rl}
\tilde{T}_{23}^j&\leq4^{1-(H+\alpha)}\|\Phi\|_{\mathcal{L}_2^{0,r+2(H+\alpha)}}^2c_H^2\left(H-\frac12\right)^2
\int_s^t\left(\int_u^t\left(r-u\right)^{2H-\frac32+\alpha}\left(\frac{r}{u}\right)^{H-\frac12}dr\right)^2du\\
&\leq\frac{4^{\frac12-(H+\alpha)}\|\Phi\|_{\mathcal{L}_2^{0,r+2(H+\alpha)}}^2}{2H+\alpha}
\left(\frac{c_H\left(H-\frac12\right)}{2H-\frac12+\alpha}\right)^2(t-s)^{4H+2\alpha}.
\end{array}
\end{equation*}
Note that when $H>\frac12$ the kernel is null on the diagonal, its
derivative has constant sign and it is integrable thus we céan
obtain without the assumption on $\alpha$
\begin{equation*}
\begin{array}{rl}
\tilde{T}_{23}^j&\leq4\int_s^t\left\|\Phi e_j \right\|_{\xHn^{1}}^2\left(\int_u^t |K(dr,u)|\right)^2du\\
&\leq4\left\|\Phi e_j
\right\|_{\xHn^{1}}^2\int_s^tK(t,u)^2du\\
&\leq4\left\|\Phi e_j
\right\|_{\xHn^{1}}^2\E\left[|\beta^H(t)-\beta^H(s)|^2\right]\\
&\leq4\left\|\Phi e_j \right\|_{\xHn^{1}}^2|t-s|^{2H}.
\end{array}
\end{equation*}
Thus $Z$ admits a modification with $H''-$H\"older continuous sample
paths with
$H'<H''<H$.\\
\indent We now explain why $Z$ has a modification which is in
$\xCn\left([0,T],\xHn^{1+2H}\right)$. Since the group is an isometry
we have
\begin{equation*}\begin{array}{rl}
\left\|Z(t)-Z(s)\right\|_{\xHn^{1+2H}}\leq&\left\|\left(U(t-s)-I\right)\sum_{j\in\xN}\int_0^T\left(K_T^*\indic_{[0,t]}(\cdot)U(-\cdot)\Phi e_j\right)(w)d\beta_j(w)\right\|_{\xHn^{1+2H}}\\
&+\left\|\tilde{T}_2(t,s)\right\|_{\xHn^{1+2H}}.
\end{array}\end{equation*}
Since the group is strongly continuous and since, from the above,
\begin{equation*}
\sum_{j\in\xN}\int_0^T\left(K_T^*\indic_{[0,t]}(\cdot)U(-\cdot)\Phi
e_j\right)(w)d\beta_j(w)
\end{equation*} belongs to $\xHn^{1+2H}$, the
first term of the right hand side goes to zero as $s$ converges to
$t$. Also, we may write
\begin{equation*}
\left\|\tilde{T}_2(t,s)\right\|_{\xHn^{1+2H}}\leq\left\|Y(t)-Y(s)\right\|_{\xHn^{1+2H}}
\end{equation*}
where $\left(Y(t)\right)_{t\in[0,T]}$, defined for $t\in[0,T]$ by
\begin{equation*}
Y(t)=\sum_{j\in\xN}\int_0^T\left(K_t^*U(-\cdot)\Phi
e_j\right)(w)d\beta_j(w),
\end{equation*}
is a Gaussian process. We again conclude, with the same bounds for
$\tilde{T}^j_{21}$ and $\tilde{T}^j_{22}$ and an upper of the order
of $(t-s)^{2H+2\alpha}$ for $\tilde{T}^j_{23}$ and using the
Kolmogorov criterion, that $Y(t)$ admits a modification with
continuous sample paths. Thus, for such a modification of $Y$, $Z$
has continuous sample paths.\\
\indent The fact that $\mu^{Z,T}$ are Gaussian measures follows from
the fact that $Z$ is defined as
\begin{equation*}
\sum_{j\in\xN}\int_0^t\left(K_T^*\indic_{[0,t]}(\cdot)U(t-\cdot)\Phi
e_j\right)(s)d\beta_j(s).
\end{equation*}
The law is Gaussian since the law of the action of an element of the
dual is a pointwise limit of Gaussian random variables; see for
example \cite{EG1}.\\
\indent It is a standard fact to prove that the process defines a
$C_{\infty}^{H,0}$ random variable, see for example \cite{EG1} for
similar arguments. We use the fact that the process takes its values
in a separable metrisable space.\hfill$\square$
\begin{rmrk}
The assumption on $\alpha$ seems too strong to have the desired
H\"older exponent. It is required only for integrability in the
upper bounds of $T_2$ and $\tilde{T}_{23}^j$. Also, the assumption
that $\Phi$ is Hilbert-Schmidt in a Sobolev space of exponent at
least $1+2H$ is only required in order that the convolution is a
$\xHn^{1+2H}$ valued process. Indeed, there is a priori no reason
that a H\"older continuous stochastic convolution gives rise to a
H\"older continuous solution to the stochastic NLS equations.
H\"older continuity of the deterministic free flow and convolution
of the nonlinearity is obtained by assuming extra space regularity
of the solution.
\end{rmrk}
\indent In the following we always consider such a modification. The
following lemma allows to characterize the RKHS of such Gaussian
measures.
\begin{lmm}\label{l2}
The covariance operator of $Z$ on $\xLtwo\left(0,T;\xLtwo\right)$ is
given for $h$ in $\xLtwo\left(0,T;\xLtwo\right)$ by
\begin{equation*}
\begin{array}{rl}
\mathcal{Q}h(t)=\sum_{j\in\xN}\int_0^T\int_0^{t\wedge u}&\left(K_T^*\indic_{[0,t]}(\cdot)U(t-\cdot)\Phi e_j\right)(s)\\
&\left(\left(K_T^*\indic_{[0,u]}(\cdot)U(u-\cdot)\Phi
e_j\right)(s),h(u)\right)_{\xLtwo}dsdu,
\end{array}\end{equation*}
when $H>\frac12$ we may write $\mathcal{Q}h(t)$ as
\begin{equation*}
c_H^2\left(H-\frac12\right)^2\beta\left(2-2H,H-\frac{1}{2}\right)\int_0^T\int_0^t\int_0^s|u-v|^{2H-2}U(t-v)\Phi\Phi^*U(u-s)h(s)dudvds.\end{equation*}
Also, for $T$ positive and $0<H'<H$, the RKHS of $\mu^{Z,T,H'}$ is
$\xim\ \mathcal{Q}^{\frac12}$ with the norm of the image structure.
It is also $\xim\ \mathcal{L}$ where $\mathcal{L}$ is defined for
$h$ in $\xLtwo\left(0,T;\xLtwo\right)$ by
\begin{equation*}
\mathcal{L}h(t)=\sum_{j\in\xN}\int_0^t\left(K_T^*\indic_{[0,t]}(\cdot)U(t-\cdot)\Phi
e_j\right)(s)(h(s),e_j)_{\xLtwo}ds.
\end{equation*}
\end{lmm}
\noindent {\bf Proof.} We may first check with the same computations
as those used in Lemma \ref{l1} that $\mathcal{L}$ is well defined
and that for $h$ in $\xLtwo(0,T;\xLtwo)$, $\mathcal{L}h$ belongs to
$\xLtwo(0,T;\xLtwo)$.
Take $h$ and $k$ in $\xLtwo(0,T;\xLtwo)$, we have
\begin{equation*}
\begin{array}{l}
\E\left[\int_0^T\left(Z(u),h(u)\right)_{\xLtwo}du\int_0^T\left(Z(t),k(t)\right)_{\xLtwo}dt \right]\\
=\sum_{j\in\xN}\E\left[\int_0^T\int_0^T\left(\int_0^T\left(K_T^*\indic_{[0,u]}(\cdot)U(u-\cdot)\Phi e_j\right)(s)d\beta_j(s),h(u) \right)_{\xLtwo}\right.\\
\quad\quad\quad\quad\quad\quad\quad\quad\left.\left(\int_0^T\left(K_T^*\indic_{[0,t]}(\cdot)U(t-\cdot)\Phi e_j\right)(v)d\beta_j(v),k(t) \right)_{\xLtwo}\right]\\
=\int_0^T\left(\mathcal{Q}h(t),k(t)\right)_{\xLtwo}dt
\end{array}
\end{equation*}
where $\mathcal{Q}$ is defined in the lemma.
The result for $H>\frac12$ is obtained with the particular form of the inner product in $\mathcal{H}$ for such values of $H$.\\
\indent Checking that for $k$ in $\xLtwo(0,T;\xLtwo)$,
\begin{equation*}
\mathcal{L}^*k(s)=\sum_{j\in\xN}\int_s^T\left(\left(K_T^*\indic_{[0,t]}
(\cdot)U(t-\cdot)\Phi e_j\right)(s),k(t)\right)_{\xLtwo}e_j dt,
\end{equation*}
we obtain that $\mathcal{Q}=\mathcal{L}\mathcal{L}^*$.\\
\indent We may thus deduce, see for example \cite{EG1}, that the
RKHS of $\mu^{Z,T,H'}$ is also $\xim\ \mathcal{L}$ with the norm of
the image structure. It is indeed the RKHS of the direct image of
$\mu^{Z,T,H'}$ on $\xLtwo(0,T;\xLtwo)$ but it is standard fact, see
for example \cite{DPZ,EG1}, that the two measures have same
RKHS.\hfill$\square$\vspace{0.3cm}

\indent When we impose (N2) we can prove as above that the
stochastic convolution $Z$ has a modification in
$\xCn\left([0,\infty);\xHn^2\right)$ embedded with the projective
limit topology letting the time interval go to infinity. Thus from
the Sobolev embeddings, for any $T$ positive and
$\left(r(p),p\right)$ an admissible pair, $Z$ belongs to
$X^{(T,p)}=\xCn\left([0,T];\xHone\right)\cap\xLn^{r(p)}\left(0,T;\xW^{1,p}\right)$.
As mentioned previsouly, this space is considered to do the fixed
point that allows to prove the local well-posedness for Kerr
nonlinearities. We may also check
\begin{lmm}\label{l2b} $Z$ defines a $\xCn\left([0,\infty);\xHn^2\right)$-random variable.
The law of its projections $\mu^{Z,T}$ on
$\xCn\left([0,T];\xHn^2\right)$ for $T$ positive is a centered
Gaussian measure whose RKHS is $\xim\
\mathcal{L}$.\end{lmm}\vspace{0.3cm}

We now deduce the following results that we will push forward to
obtain results for the solution of the SPDE.
\begin{prpstn}\label{p1}
The direct image measures for $\epsilon$ positive of
$x\mapsto\sqrt{\epsilon}x$ on $C_{\infty}^{H,0}$, respectively
$\xCn\left([0,\infty);\xHn^2\right)$, satisfy a LDP of speed
$\epsilon$ and good rate function
\begin{equation*}
I^Z(f)=\frac12\inf_{h\in\xLtwo\left(0,\infty;\xLtwo\right):\
\mathcal{L}(h)=f}\left\{\|h\|_{\xLtwo\left(0,\infty;\xLtwo\right)}^2\right\}.
\end{equation*}
\end{prpstn}
\noindent {\bf Proof.} From a general result on LDP for Gaussian
measures on Banach spaces, see \cite{DS1}, and the above lemma, we
know that for $T$ positive and $0<H'<H$, the direct images of
$\mu^{Z,T,H'}$ by the mapping $x\mapsto\sqrt{\epsilon}x$ satisfy a
LDP of speed $\epsilon$ and good rate function
\begin{equation*}
I^{Z,T,H'}(f)=\frac12\inf\left\{ \|h\|_{\xim\ \mathcal{L}}^2:\
f=\mathcal{L}h\right\}
\end{equation*}
with the convention that $\inf \emptyset=\infty$. We conclude
letting $T$ go to infinity and $H'$ to $H$ using Dawson-Gartner's
theorem for projective limits, see for example \cite{DZ}, and
Lebesgue's dominated convergence theorem. The same is true under
(N2) when we work in
$\xCn\left([0,\infty);\xHn^2\right)$.\hfill$\square$
\begin{prpstn}\label{p2}
Under (N1) the support of the measure $\mu^{Z}$ is given by
\begin{equation*}
{\rm supp}\ \mu^{Z}= \overline{\xim\
\mathcal{L}}^{C_{\infty}^{H,0}},
\end{equation*}
under (N2) the same result holds replacing $C_{\infty}^{H,0}$ by
$\xCn\left([0,\infty);\xHn^2\right)$.
\end{prpstn}
\noindent {\bf Proof.} Let us give the argument when we have the
assumption (N1), the argument under (N2) is the same.\\
From the characterization of the RKHS of the centered Gaussian
measure $\mu^{Z,T,H'}$ for $T$ positive and $0<H'<H$ and Theorem
(IX,2;1) in \cite{BC}, we obtain that the support of the measure
$\mu^{Z,T,H'}$ is such that
\begin{equation*}
{\rm supp}\ \mu^{Z,T,H'}= \overline{\xim\ \mathcal{L}}^{C_T^{H',0}}.
\end{equation*}
From the definition of the image measure we have that
\begin{equation*}
\mu^Z\left(p_{T,H'}^{-1}\left(\overline{\xim\
\mathcal{L}}^{C_T^{H',0}}\right)\right)=\mu^{Z,T,H'}\left(\overline{\xim\
\mathcal{L}}^{C_T^{H',0}}\right)=1,
\end{equation*}
where $p_{T,H'}$ denotes the projection of $C_{\infty}^{H,0}$ into
$C_T^{H',0}$. It follows that
\begin{equation*}{\rm supp}\
\mu^Z\subset\bigcap_{T}p_{T,H'}^{-1}\left(\overline{\xim\
\mathcal{L}}^{C_T^{H',0}}\right)=\overline{\xim\
\mathcal{L}}^{C_{\infty}^{H,0}}.
\end{equation*}
It then suffices to show that $\xim\ \mathcal{L}\subset {\rm supp}\
\mu^Z$. Suppose that $x\notin {\rm supp}\ \mu^Z$, then there exists
a neighborhood $V$ of $x$ in $C_{\infty}^{H,0}$ which is a
neighborhood of $x$ in $C_T^{H',0}$ for $T$ large and $H'$
sufficiently close to $H$ such that $\mu^Z(V)=0$. Since the support
of $\mu^{Z,T,H'}$ is the closure of $\xim\ \mathcal{L}$ for the
topology of $C_T^{H',0}$, $V\cap\xim\ \mathcal{L}=\emptyset$ and
$x\notin \xim\ \mathcal{L}$.\hfill$\square$

\section{Local well-posedness of the Cauchy problem}
We consider the Cauchy problem
\begin{equation}\label{Cauchy}
\left\{\begin{array}{l} i\xdif u=\left(\Delta u +f(u)\right)\xdif t+dW^H\\
u(0)=u_0.
\end{array}\right.
\end{equation}
We consider two cases. In the first case we assume (N1),
$u_0\in\xHn^{1+2H}$ and\vspace{0.2cm}

\noindent {\bf Assumption (NL)}
\begin{equation*}
\begin{array}{rl}
(i)&f\ {\it is\ Lipschitz\ on\ the\ bounded\ sets\ of}\
\xHn^{1+2H}\\
(ii)&f(0)=0.
\end{array}
\end{equation*}
In the second case we assume (N2), $u_0\in\xHone$ and $f$ is a Kerr
nonlinearity.\vspace{0.3cm}

\indent We first recall the following important fact. Let us denote
by $v^{u_0}(z)$ the solution of
\begin{equation}\label{Cauchy2}
\left\{\begin{array}{l} i\frac{\xdif v}{\xdif t}=\Delta v +f(v-iz)\\
u(0)=u_0.
\end{array}\right.
\end{equation}

where $z$ is a function of $C_{\infty}^{H,0}$ (repectively
$\xCn\left([0,\infty),\xHn^2\right)$) and define $\mathcal{G}^{u_0}$
the mapping
\begin{equation*}\mathcal{G}^{u_0}:\ z\mapsto
v^{u_0}(z)-iz.
\end{equation*}
Then we may check that the solution $u^{\epsilon,u_0}$ of
\eqref{Cauchy} is such that
$u^{\epsilon,u_0}=\mathcal{G}^{u_0}(\sqrt{\epsilon}Z)$ where $Z$ is
the stochastic convolution.\\
\indent We may now check with a fixed point argument the following
result.
\begin{thrm}\label{t1}
Assume that the initial datum $u_0$ is $\mathcal{F}_0$ measurable
and belongs to $\xHn^{1+2H}$ (repsectively $\xHone$); then there
exists a unique solution to \eqref{e2} with continuous $\xHn^{1+2H}$
(respectively $\xHone$) valued paths. The solution is defined on a
random interval $[0,\tau^*(u_0,\omega))$ where $\tau^*(u_0,\omega)$
is either $\infty$ or a finite blow-up time.
\end{thrm}
In the next two section we state sample paths LDPs and support
theorems. We start with the first set of assumptions and state a
result in a space of H\"older continuous sample paths with any value
of the Hurst parameter. In the last section we consider the case of
Kerr nonlinearities and restrict ourselves to the case where
$H>\frac12$.

\section{The case of a nonlinearity satisfying (NL)}\label{s4}
According to (NL) solutions may blow up in finite time. We shall
proceed as in \cite{EG1} to define proper path spaces where we can
state the LDP and support result; see the reference for more
details. However, we consider here a space where paths are
$H'-$H\"older continuous with values in $\xHn^1$ on compact time
intervals before the blow-up time where $0<H'<H$. We add a point
$\Delta$ to the space $\xHn^{1+2H}$ and embed the space with the
topology such that its open sets are the open sets of $\xHn^{1+2H}$
and the complement in $\xHn^{1+2H}\cup\{\Delta\}$ of the closed
bounded sets of $\xHn^{1+2H}$. The set
$\xCn([0,\infty);\xHn^{1+2H}\cup\{\Delta\})$ is then well defined.
We denote the blow-up time of $f$ in
$\xCn([0,\infty);\xHn^{1+2H}\cup\{\Delta\})$ by
$\mathcal{T}(f)=\inf\{t\in[0,\infty):\ f(t)=\Delta\}$, with the
convention that $\inf\emptyset=\infty$.\\
\indent We also define the following spaces
\begin{equation*}
C_T^{H'}=\xCn\left([0,T];\xHn^{1+2H}\right)\cap\xCn^{H'}\left([0,T];\xHn^{1}\right)
\end{equation*}
and
\begin{equation*}
C_T^{H',0}=\xCn\left([0,T];\xHn^{1+2H}\right)\cap\xCn^{H',0}\left([0,T];\xHn^{1}\right).
\end{equation*}
When equipped with the norm which is the supremum of the norms of
the two Banach spaces intersected they are Banach spaces. The latter
space is separable.\\
\indent For measurability issue, we define
\begin{equation*}
C_{\infty}^{H,0}=\bigcap_{T>0,0<H'<H}C_T^{H',0}
\end{equation*}
equipped with the projective limit topology. It is a separable
metrisable space. We also define
\begin{equation*}\begin{array}{rl}
\mathcal{E}^{H}\left(\xHn^{1}\right)=&\left\{f\in\xCn([0,\infty);\xHn^{1+2H}\cup\{\Delta\}):\
f(t_0)=\Delta\Rightarrow\forall t\geq t_0,\ f(t)=\Delta ;\right.\\
&\left.\ \ \forall T<\mathcal{T}(f),\ \forall\ 0<H'<H,\ f\in
\xCn^{H'}\left([0,T];\xHn^{1}\right)\right\}.
\end{array}\end{equation*}
Here $\Delta$ acts as a cemetary. It is endowed with the topology
defined by the neighborhood basis
\begin{equation*}
V_{T,R,H'}(\varphi_1)=\left\{\varphi\in\mathcal{E}^{H}\left(\xHn^{1}\right):\
\mathcal{T}(\varphi)> T,\ \|\varphi_1-\varphi\|_{C_T^{H'}}\leq
R\right\},
\end{equation*}
of $\varphi_1$ in $\mathcal{E}^{H}\left(\xHn^1\right)$ given
$T<\mathcal{T}(\varphi_1)$ and $R$ positive. The space is also a
Hausdorff topological space and thus we may consider applying the
Varadhan contraction principle.\\
\indent In order to push forward the results of section \ref{s3} we
use the following result.
\begin{lmm}\label{l3}
The mapping
\begin{equation*}
\begin{array}{rrl}
C_{\infty}^{H,0}&\rightarrow&\mathcal{E}^{H}\left(\xHn^{1}\right)\\
z&\mapsto&\mathcal{G}^{u_0}(z)
\end{array}
\end{equation*}
is continuous.
\end{lmm}
\noindent {\bf Proof.} This could be done by revisiting the fixed
point argument, this time in $C_{T^*}^{H'}$ for $T^*$ small enough
depending on the norm of the initial data and $z$ in $C_{T}^{H'}$
for some fixed $T$ and some $H'<H$ fixed. Though with different
norms, the remaining of the argument allowing to prove the
continuity of $v^{u_0}(z)$ with respect to $z$, detailed in
\cite{dBD1}, holds. In the computations we use \eqref{PG} in order
to treat the H\"older norms. \hfill$\square$\vspace{0.3cm}

Note that Lemma \ref{l1} and \ref{l3} give that $u^{1,u_0}$ defines
a $\mathcal{E}^{H}\left(\xHn^{1}\right)$ random
variable.\vspace{0.3cm}

\indent Let us now study large deviations for the laws
$\mu^{u^{\epsilon,u_0}}$ on $\mathcal{E}^{H}(\xHn^{1})$ of the mild
solutions $u^{\epsilon,u_0}$ of
\begin{equation}\label{e6}
\left\{\begin{array}{l} i \xdif u -(\Delta u +
f(u))\xdif t=\sqrt{\epsilon}\xdif W^H,\\
u(0)=u_0\ \in\xHn^{1+2H}.
\end{array}\right.
\end{equation}
We may now deduce from Lemma \ref{l2} and \ref{l3}, the fact that
$\left(\mathcal{G}^{u_0}\circ\mathcal{L}\right)(\cdot)=\mathbf{S}(u_0,\cdot)$,
and the Varadhan contraction principle the following theorem.
\begin{thrm}\label{t2}
The laws $\mu^{u^{\epsilon,u_0}}$ on
$\mathcal{E}^{H}\left(\xHn^1\right)$ satisfy a LDP of speed
$\epsilon$ and good rate function
\begin{equation*}
I^{u_0}(w)=\frac12\inf_{h\in\xLtwo\left(0,\infty;\xLtwo\right):\
\mathbf{S}(u_0,h)=w}
\left\{\|h\|_{\xLtwo\left(0,\infty;\xLtwo\right)}^2\right\},
\end{equation*}
where $\mathbf{S}(u_0,h)$ denotes the mild solution in
$\mathcal{E}^{H}\left(\xHn^1\right)$ of the following control
problem
\begin{equation}\label{control}
\left\{\begin{array}{l}
i \frac{\partial u}{\partial t} -(\Delta u + f(u)) =\Phi \mbox{\.{Kh}},\\
u(0)=u_0\ \in\xHn^{1+2H},\ h\in \xLtwo\left(0,\infty;\xLtwo\right);
\end{array}\right.
\end{equation}
it is called the skeleton. Only the integral, or the integral in the mild formulation, of the right hand side
is defined; it is by means of the duality relation.
\end{thrm}
\begin{rmrk}\label{rk2}
We could also prove a uniform LDP as for example in \cite{EG2}.
\end{rmrk}
The characterization of the support follows with the same arguments
as in \cite{EG1}. We recall the proof for the sake of completeness.
\begin{thrm}\label{t3}
The support of the law $\mu^{u^{1,u_0}}$ on
$\mathcal{E}^{H}\left(\xHn^{1}\right)$ is given by
\begin{equation*}
{\rm supp}\ \mu^{u^{1,u_0}}=\overline{\xim\
\mathbf{S}}^{\mathcal{E}^{H}\left(\xHn^{1}\right)}.
\end{equation*}
\end{thrm}
\noindent {\bf Proof.} We use the continuity of $\mathcal{G}$.
Indeed, since $\mathcal{G}^{u_0}(\xim\
\mathcal{L})\subset\overline{\mathcal{G}^{u_0}(\xim\
\mathcal{L})}^{\mathcal{E}^{H}\left(\xHn^{1}\right)}$, $\xim\
\mathcal{L}\subset\left(\mathcal{G}^{u_0}\right)^{-1}\left(\overline{\mathcal{G}^{u_0}(\xim\
\mathcal{L})}^{\mathcal{E}^{H}\left(\xHn^{1}\right)}\right)$.
Because $\mathcal{G}^{u_0}$ is continuous, the right hand side is a
closed set of $C_{\infty}^{H,0}$ and from Proposition \ref{p2},
\begin{equation*}{\rm supp}\
\mu^Z\subset\left(\mathcal{G}^{u_0}\right)^{-1}\left(\overline{\xim\
\left(
\mathcal{G}^{u_0}\circ\mathcal{L}\right)}^{\mathcal{E}^{H}\left(\xHn^{1}\right)}\right),
\end{equation*}
and
\begin{equation*}\mu^Z\left(\left(\mathcal{G}^{u_0}\right)^{-1}\left(\overline{\xim\
\mathbf{S}(u_0)}^{\mathcal{E}^{H}\left(\xHn^{1}\right)}\right)\right)=1,
\end{equation*}
thus
\begin{equation*}{\rm supp}\ \mu^u\subset\overline{\xim\
\mathbf{S}(u_0)}^{\mathcal{E}^{H}\left(\xHn^{1}\right)}.
\end{equation*} Suppose that
$x\notin {\rm supp}\ \mu^{u^{1,u_0}}$, there exists a neighborhood
$V$ of $x$ in $\mathcal{E}^{H}\left(\xHn^{1}\right)$ such that
$\mu^{u^{1,u_0}}(V)=\mu^Z\left(\left(\mathcal{G}^{u_0}\right)^{-1}(V)\right)=0$,
consequently $\left(\mathcal{G}^{u_0}\right)^{-1}(V)\bigcap \xim\
\mathcal{L}$ is empty and $x\notin \xim\ \mathbf{S}(u_0)$. This
gives the reverse inclusion.\hfill $\square$

\section{The case of Kerr nonlinearities}\label{s5}
In this section we consider Kerr nonlinearities when
$d\geq2$ and $\sigma<\frac{2}{d-2}$.\\
\indent This time, we will not state a result in a space of H\"older
continuous functions with values in $\xHone$. We would need that the
convolution which involves the nonlinearity is H\"older continuous.
Thus, in order to use \eqref{PG}, we would have to compute the
Sobolev norm of the nonlinearity in some space $\xHn^{1+2\gamma}$
where $\gamma$ is positive.
\begin{rmrk}
In the case where $H<\frac12$, we could however state a weaker
result than in the previous section imposing that $u_0\in\xHone$ and
$\Phi\in\mathcal{L}_2^{0,2+\alpha}$. The corresponding fixed point
could be conducted in
$\xCn^{H'}\left([0,T];\xHn^{1-2H}\right)\cap\xCn\left([0,T];\xHone\right)\cap\xLn^{r(p)}\left(0,T;\xW^{1,p}\right)$
where $\left(r(p),p\right)$ is an admissible pair and $0<H'<H$ and
uses the Strichartz inequalities. Indeed, from the Sobolev
embeddings, the stochastic convolution has a modification in
$\xCn\left([0,T];\xHn^2\right)\cap\xCn^{H'}\left([0,T];\xHn^{2-2H}\right)$
and thus belong to the desired space. \end{rmrk} Let us return to
the case where $H>\frac12$. Since under (N2) we know that the
stochastic convolution $Z$ has a modification in $X^{(T,p)}$, we can
directly use the continuity of the solution with respect to the
stochastic convolution of \cite{EG1} and repeat the remining of the
arguments. Thus, for initial data in $\xHone$, we may state a LDP
and support result in the space $\mathcal{E}_{\infty}$ defined as
\begin{equation*}\begin{array}{l}
\mathcal{E}_{\infty}=\left\{f\in\xCn([0,\infty);\xHone\cup\{\Delta\}):\
f(t_0)=\Delta\Rightarrow\forall t\geq t_0,\ f(t)=\Delta ;\right.\\
\hspace{1.2cm} \left.\forall T<\mathcal{T}(f),\ \forall\ p\in
\left[2,\frac{2d}{d-2}\right), f\in
\xL^{r(p)}\left(0,T;\xW^{1,p}\right)\right\}.
\end{array}\end{equation*}
When $d=2$ or $d=1$ we write $p\in[2,\infty)$. The space is embedded
with the topology defined by the neighborhood basis
\begin{equation*}
W_{T,p,R}(\varphi_1)=\left\{\varphi\in\mathcal{E}_{\infty}:
\mathcal{T}(\varphi)\geq T,\ \|\varphi_1-\varphi\|_{X^{(T,p)}}\leq
R\right\},
\end{equation*}
for $\varphi_1$ in $\mathcal{E}_{\infty}$.
\begin{thrm}\label{t4}
The laws $\mu^{u^{\epsilon,u_0}}$ on $\mathcal{E}_{\infty}$ satisfy
a LDP of speed $\epsilon$ and good rate function
\begin{equation*}
I^{u_0}(w)=\frac12\inf_{h\in\xLtwo\left(0,\infty;\xLtwo\right):\
\mathbf{S}(u_0,h)=w}
\left\{\|h\|_{\xLtwo\left(0,\infty;\xLtwo\right)}^2\right\},
\end{equation*}
where $\mathbf{S}(u_0,h)$ is the mild solution of
\begin{equation}\label{control2}
\left\{\begin{array}{l}
i \frac{\partial u}{\partial t} -(\Delta u +\lambda|u|^{2\sigma}u) =\Phi \mbox{\.{Kh}},\\
u(0)=u_0\ \in\xHn^{1},\ h\in \xLtwo\left(0,\infty;\xLtwo\right);
\end{array}\right.
\end{equation}
\end{thrm}
\begin{thrm}\label{t5}
The support of the law $\mu^{u^{1,u_0}}$ on $\mathcal{E}_{\infty}$
is given by
\begin{equation*}
{\rm supp}\ \mu^{u^{1,u_0}}=\overline{\xim\
\mathbf{S}}^{\mathcal{E}_{\infty}}.
\end{equation*}
\end{thrm}

\end{document}